\documentclass[12pt]{article}
\usepackage{amsfonts}
\usepackage{amsmath,amssymb,amsfonts}
\usepackage{color}
\usepackage{hyperref}
\usepackage[latin1]{inputenc}
\usepackage{amsthm}
\usepackage{graphicx}
\usepackage{ulem}

\def\squarebox#1{\hbox to #1{\hfill\vbox to #1{\vfill}}}
\newcommand{\qedem}{\hspace*{\fill}
\vbox{\hrule\hbox{\vrule\squarebox{.667em}\vrule}\hrule}\smallskip}
\newtheorem{teo}{Theorem}[section]
\newtheorem{prop}[teo]{Proposition}

\newtheorem{coro}[teo]{Corollary}

\newtheorem{lema}[teo]{Lemma}

\newtheorem{rem}[teo]{Remark}
\newenvironment{profe}{\noindent {\bf Proof:}}{\hfill $\qedem $ \newline}

\begin{document}

\title{Semigroups and Controllability of Invariant Control Systems on $%
\mathrm{Sl}\left( n,\mathbb{H}\right) $}
\author{Bruno A. Rodrigues\\Departamento de Matem\'{a}tica, Universidade Estadual de Maring\'{a}, Brazil \and Luiz A. B. San Martin\\Departamento de Matem\'{a}tica, Universidade de Campinas, Brazil \and  Alexandre J. Santana\\Departamento de Matem\'{a}tica, Universidade Estadual de Maring\'{a}, Brazil}

\maketitle

\textbf{Abstract: }
Let $\mathrm{Sl}\left( n,\mathbb{H}\right) $ be the Lie group of $n\times
n $ quaternionic matrices $g$ with $\left\vert \det g\right\vert =1$. We prove that a 
subsemigroup $S \subset \mathrm{Sl}\left( n,\mathbb{H}\right)$ with nonempty interior is equal to $\mathrm{Sl}\left( n,\mathbb{H}\right)$ if $S$ contains a subgroup isomorphic to $\mathrm{Sl}\left( 2,\mathbb{H}\right)$. As application we
give
sufficient conditions on $A,B\in \mathfrak{sl}\left( n,%
\mathbb{H}\right) $ to ensuring that the invariant control system $\dot{g}%
=Ag+uBg$ is controllable on $\mathrm{Sl}\left( n,\mathbb{H}\right) $. We prove also that these
conditions are generic in the sense that we obtain an open and dense set of controllable  pairs $\left( A,B\right)
\in \mathfrak{sl}\left( n,\mathbb{H}\right) ^{2}$.


\textbf{MSC\ 2010.} 93B05, 22E20, 22F30.

\section{Introduction and main results}

Let $\mathrm{Sl}\left( n,\mathbb{H}\right) $ be the simple Lie group of
quaternionic $n\times n$ matrices $g$ with $\left\vert \det g\right\vert =1$. In this paper we give conditions ensuring that a subsemigroup $S\subset 
\mathrm{Sl}\left( n,\mathbb{H}\right) $ is in fact the whole group. These
conditions are applied to the controllability problem of an invariant
control system%
\begin{equation}
\dot{g}=Ag+uBg  \label{forsis}
\end{equation}%
on $\mathrm{Sl}\left( n,\mathbb{H}\right) $. In such a system, $%
A $ and $B$ are quaternionic $n \times n$ matrices  having trace with zero real part, i.e.,
$A,B$ are elements of the Lie algebra $\mathfrak{sl}\left( n,\mathbb{H}
\right) $ of $\mathrm{Sl}\left( n,\mathbb{H}\right) $.  This is
a real simple Lie algebra that complexifies to a complex Lie algebra
isomorphic to $\mathfrak{sl}\left( 2n,\mathbb{C}\right) $.

Our approach to prove that a subsemigroup $S$ equals $\mathrm{Sl}\left( n,%
\mathbb{H}\right) $ follows the same lines of the results proved by Dos Santos and San Martin in \cite%
{smals} and \cite{smals1}. These papers work in a general noncompact
semi-simple Lie group $G$ and develop a method based on the topology of flag
manifolds of $G$. That topological method permits to show that a
subsemigroup $S$ with $\mathrm{int}S\neq \emptyset $ must be the group $G$
provided it contains certain classes of subgroups given by a root of the Lie algebra of $G$.

For the group $\mathrm{Sl}\left( n,\mathbb{H}\right) $,   our main result proves that  
 $
S=\mathrm{Sl}\left( n,\mathbb{H}\right) $ if $\mathrm{int}S\neq \emptyset $
and $S$ contains the subgroup $\langle \exp \mathfrak{g}_{\pm \alpha
}\rangle $ generated by the root spaces $\mathfrak{g}_{\pm \alpha }$ of a
root $\alpha $. Explicitly, for a pair $\left( r,s\right) $, $1\leq r<s\leq n$
let $\mathrm{Sl}\left( 2,\mathbb{H}\right) _{r,s}$ be the subgroup of $%
\mathrm{Sl}\left( n,\mathbb{H}\right) $, isomorphic to $\mathrm{Sl}\left( 2,%
\mathbb{H}\right) $, of the matrices in the space $\mathrm{span}%
\{e_{r},e_{s}\}$ where $\{e_{1},\ldots ,e_{n}\}$ is the standard basis of $%
\mathbb{H}^{n}$, plugged into the $n\times n$ matrices. That is, $\mathrm{Sl}%
\left( 2,\mathbb{H}\right) _{r,s}=\langle \exp \mathfrak{sl}\left( 2,\mathbb{%
	H}\right) _{r,s}\rangle $ where the Lie algebra $\mathfrak{sl}\left( 2,%
\mathbb{H}\right) _{r,s}$ is given by the matrices in $\mathfrak{sl}\left( n,%
\mathbb{H}\right) $ having nonzero entries only at the positions $\left(
a,b\right) $ with $a,b\in \{r,s\}$. 

Then we apply this theorem to the controllability of (\ref{forsis}). We
follow an idea that goes back to the papers by Jurdjevic and Kupka \cite{JK}
and \cite{JK1}. More precisely, we find conditions on the matrices $A$ and $B$ ensuring
that $\mathrm{Sl}\left( 2,\mathbb{H}\right) _{1,n}$ is contained in the
control semigroup $S_{A,B}$, that is, the subsemigroup of $\mathrm{Sl}\left(
n,\mathbb{H}\right) $ generated by the $1$-parameter semigroups $e^{t\left(
A+uB\right) }$, $t\geq 0$, $u\in \mathbb{R}$. Hence if we assume also the
Lie algebra rank condition (so that $\mathrm{int}S_{A,B}\neq \emptyset $)
then Theorem \ref{teoTopology} applies to get $S_{A,B}=\mathrm{Sl}\left( n,%
\mathbb{H}\right) $, which means that (\ref{forsis}) is controllable. 

In this
context we prove that controllability for invariant systems on $\mathrm{Sl}%
\left( n,\mathbb{H}\right) $ is a generic property in the sense that there is an open and 
 dense set $C\subset \mathfrak{sl}%
\left( n,\mathbb{H}\right) ^{2}$ such that the control system $\dot{g}%
=A\left( g\right) +uB\left( g\right) $ with unrestricted controls ($u\in 
\mathbb{R}$) is controllable for all pairs $\left( A,B\right) \in C$. 
To get the set $C\subset \mathfrak{sl}\left( n,\mathbb{H}\right) ^{2}$  we first write a sufficient condition for controllability by
taking $B$ to be diagonal (see Theorem \ref{teoregularreal} below).
Afterwards we check that in $\mathfrak{sl}\left( n,\mathbb{H}\right) $ there
is just one conjugacy class of Cartan subalgebras. This allows to build the
dense set $C$ as a set of conjugates of a controllable pair $\left(
A,B\right) $ with $B$ diagonal (see Section \ref{secgeneric} below).

About the structure of this paper, we first present the fundamental concepts for this paper. In the third section, we develop and prove our
main result  which gives necessary conditions for a subsemigroup of $\mathrm{Sl}%
\left( n,\mathbb{H}\right) $ to be equal to $\mathrm{Sl}%
\left( n,\mathbb{H}\right) $. In the fourth section we apply the previous results to prove
 a theorem that provides sufficient conditions for controllability of invariant systems on  $\mathrm{Sl}%
 \left( n,\mathbb{H}\right) $. Finally, in the last section we show that controllability, for the above system, is a generic property.

\section{Background}

In this section we establish  some necessary notations, concepts and results.

In a matrix Lie algebra a  Cartan decomposition
$\mathfrak{g}=\mathfrak{k}\oplus \mathfrak{s}$ is given by skew symmetric  and symmetric (or hermitian) matrices. Hence the  natural Cartan decomposition of $\mathfrak{sl}\left( n,\mathbb{H}\right) $ is 
\begin{equation*}
\mathfrak{k}=\{X\in M_{n\times n}\left( \mathbb{H}\right) :X=-\overline{X}%
^{T}\}\qquad \mathfrak{s}=\{X\in M_{n\times n}\left( \mathbb{H}\right) :X=%
\overline{X}^{T}\}
\end{equation*}%
where $\overline{\cdot }$ is a quaternionic conjugation. The algebra $\mathfrak{k}
$ of the quaternionic skew Hermitian  matrices is denoted by $\mathfrak{k}=%
\mathfrak{sp}\left( n\right) $ and is the compact real form of $C_{n}=\mathfrak{%
	sp}\left( n,\mathbb{C}\right) $.

The maximal abelian subalgebra $\mathfrak{a}%
\subset \mathfrak{s}$ is given by the diagonal matrices  $%
\Lambda =\mathrm{diag}\{a_{1},\ldots ,a_{n}\}$ with $a_{i}\in \mathbb{R}$ and $%
\mathrm{tr}\Lambda =0$. The roots of $\mathfrak{a}$  are the following linear functionals  
\begin{equation*}
\alpha _{rs}\left( \Lambda \right) =\left( \lambda _{r}-\lambda _{s}\right)
\left( \Lambda \right) =a_{r}-a_{s}\qquad r\neq s.
\end{equation*}%
The vector space  $\mathfrak{g}_{\alpha _{rs}}$ corresponding to the root  $\alpha _{rs}$ is given
by the quaternionic matrices with non zero       entries only in the position  $rs$. Then all roots have multiplicity $4$.
The set of simple system of roots is given by  
\begin{eqnarray*}
	\Sigma &=&\{\lambda _{1}-\lambda _{2},\ldots ,\lambda _{n-1}-\lambda _{n}\}
	\\
	&=&\{\alpha _{1},\ldots ,\alpha _{n-1}\} .
\end{eqnarray*}%

With this choice, the set of positive roots is given by $\alpha _{rs}$ with $r<s$. Hence an Iwasawa decomposition is  
\begin{equation*}
\mathfrak{sl}\left( n,\mathbb{H}\right) =\mathfrak{sp}\left( 1\right) \oplus 
\mathfrak{a}\oplus \mathfrak{n}^{+}
\end{equation*}%
where $\mathfrak{n}^{+}$ is the Lie algebra of upper triangular quaternionic  $n\times n$ matrices  with zero entries in the diagonal.

Now we recall an important result for our goals, which is related to the flag type of a semigroup in the special case of  $\mathrm{Sl}\left( n,\mathbb{H}\right) $ (for the
general context see San Martin and Tonelli \cite{SMT}, San Martin \cite{smmax}, Dos Santos and San Martin \cite{smals1} and references
therein). For $d=1,\ldots ,n-1$, we denote by $\mathrm{Gr}_{d}\left( \mathbb{H%
}\right) $ the Grassmannian of $d$-dimensional quaternionic subspaces of $%
\mathbb{H}^{n}$. The group $\mathrm{Sl}\left( n,\mathbb{H}\right) $ acts
transitively on each $\mathrm{Gr}_{d}\left( \mathbb{H}\right) $. The compact
subgroup $\mathrm{Sp}\left( n\right) \subset \mathrm{Sl}\left( n,\mathbb{H}%
\right) $ also acts transitively on $\mathrm{Gr}_{d}\left( \mathbb{H}\right) 
$.

\begin{teo}
	\label{teoflagtype}Let $S\subset \mathrm{Sl}\left( n,\mathbb{H}\right) $ be
	a proper subsemigroup with $\mathrm{int}S\neq \emptyset $. Then there are $%
	d\in \{1,\ldots ,n-1\}$ and a subset $C_{d}\subset \mathrm{Gr}_{d}\left( 
	\mathbb{H}\right) $ satisfying
	
	\begin{enumerate}
		\item $C_{d}$ is closed, has nonempty interior and is invariant by the action of 
		$S$. ($C_{d}$ is the unique invariant control set of $S$ in $\mathrm{Gr}%
		_{d}\left( \mathbb{H}\right) $).
		
		\item $C_{d}$ is contractible in $\mathrm{Gr}_{d}\left( \mathbb{H}\right) $
		in the sense that there exists $H\in \mathfrak{sl}\left( n,\mathbb{H}\right) 
		$ such that $e^{tH}C_{d}$ shrinks to a point as $t\rightarrow +\infty $.
	\end{enumerate}
\end{teo}

In the context of the above theorem, the Grassmannian  $\mathrm{Gr}_{d}\left( \mathbb{H}\right) $ is called the flag type of the semigroup $S$.

\section{Transitivity of a subsemigroup of $\mathrm{Sl}\left( n,\mathbb{H}\right)$}

In this section, following the same construction and notation in the introduction, we prove our main result that gives sufficient condition for a semigroup $S$ to be equal to $\mathrm{Sl}\left( n,%
\mathbb{H}\right) $. The proof of this result is based on the existence of a flag type of 
a proper semigroup $S$ with $\mathrm{int}S\neq
\emptyset $ and it follows the same pattern as the proof
of the results in \cite{smals} and \cite{smals1}. By Theorem \ref{teoflagtype}, we get that $S=\mathrm{Sl}\left( n,\mathbb{H}\right) $ if we
can prove that $S$ does not leave invariant contractible subsets in the
Grassmannians $\mathrm{Gr}_{d}\left( \mathbb{H}\right) $, $d=1,\ldots ,n-1$. 

Now we state the main theorem of this paper.

\begin{teo}
	\label{teoTopology}Let $S\subset \mathrm{Sl}\left( n,\mathbb{H}\right) $ be
	a subsemigroup with $\mathrm{int}S\neq \emptyset $ and suppose that $\mathrm{%
		Sl}\left( 2,\mathbb{H}\right) _{r,s}\subset S$ for some pair of indices $%
	\left( r,s\right) $, $1\leq r<s\leq n$. Then $S=\mathrm{Sl}\left( n,\mathbb{H%
	}\right) $.
\end{teo}

Before to prove this theorem, we need some remarks and lemmas. 
Taking into account the assumption, of the above theorem,  that $%
\mathrm{Sl}\left( 2,\mathbb{H}\right) _{r,s}\subset S$ we consider
separately the case where $\left( r,s\right) =\left( 1,n\right) $, that is, $%
\mathrm{Sl}\left( 2,\mathbb{H}\right) _{1,n}\subset S$. Then the strategy
is to prove that an $S$-invariant subset $C_{d}\subset \mathrm{Gr}_{d}\left( 
\mathbb{H}\right) $, which is closed and has nonempty interior, contains an
orbit of $\mathrm{Sl}\left( 2,\mathbb{H}\right) _{1,n}$ that is not
contractible to a point in $\mathrm{Gr}_{d}\left( \mathbb{H}\right) $. In
the next lemma we describe the noncontractible $\mathrm{Sl}\left( 2,\mathbb{H%
}\right) _{1,n}$-orbits in the Grassmannians that are proved to be contained
in the invariant control sets $C_{d}$. They are $4$-dimensional spheres (cf.
Lemma 3.7 of \cite{smals}).

\begin{lema}
\label{lemEsse4}For $d=1,n-1$ let $V_{d}$ be the subspace of $\mathbb{H}^{n}$
spanned by the first $d$ basic vectors,%
\begin{equation*}
V_{d}=\{\left( q_{1},\ldots ,q_{d},0,\ldots ,0\right) :q_{r}\in \mathbb{H}\}.
\end{equation*}%
Then the $\mathrm{Sl}\left( 2,\mathbb{H}\right) _{1,n}$-orbit in $\mathrm{Gr}%
_{d}\left( \mathbb{H}\right) $ through $V_{d}$ is diffeomorphic to $S^{4}$.
\end{lema}

\begin{profe}
The orbit is diffeomorphic to the coset space $\mathrm{Sl}\left( 2,\mathbb{H}%
\right) _{1,n}/P$ where $P=\{g:gV_{d}=V_{d}\}$ is the isotropy subgroup. By
a direct check one sees that $P$ is the subgroup of matrices in $\mathrm{Sl}%
\left( 2,\mathbb{H}\right) _{1,n}$ that are upper triangular. This is a
parabolic subgroup of $\mathrm{Sl}\left( 2,\mathbb{H}\right) _{1,n}$ hence $%
\mathrm{Sl}\left( 2,\mathbb{H}\right) _{1,n}/P$ is a flag manifold of $%
\mathrm{Sl}\left( 2,\mathbb{H}\right) _{1,n}$. Now $\mathrm{Sl}\left( 2,%
\mathbb{H}\right) $ is a real rank $1$ group so that it has just one flag
manifold which is diffeomorphic to a sphere. The dimension of the sphere
equals the codimension of $P$ which is $4$. Therefore $\mathrm{Sl}\left( 2,%
\mathbb{H}\right) _{1,n}/P$ as well as the orbit through $V_{d}$ is a sphere 
$S^{4}$.
\end{profe}

The next step is to check that for any $d=1,\ldots ,n-1$ the orbit $$\mathrm{%
Sl}\left( 2,\mathbb{H}\right) _{1,n}V_{d} \approx S^{4}$$ is not contractible
in $\mathrm{Gr}_{d}\left( \mathbb{H}\right) $, that is, is not homotopic to
a point. In other words we are required to prove that the $4$-sphere $%
\mathrm{Sl}\left( 2,\mathbb{H}\right) _{1,n}V_{d}$ is not a representative
of the identity of the homotopy group $\pi _{4}\left( \mathrm{Gr}_{d}\left( 
\mathbb{H}\right) \right) $. To this purpose we recall the cellular
decomposition of $\mathrm{Gr}_{d}\left( \mathbb{H}\right) $ given in Rabelo and San Martin \cite%
{papercell}. From that decomposition the homology $H_{\ast }\left( \mathrm{Gr%
}_{d}\left( \mathbb{H}\right) \right) $ of a Grassmannian $\mathrm{Gr}%
_{d}\left( \mathbb{H}\right) $ is freely generated by the Schubert cells and 
$H_{r}\left( \mathrm{Gr}_{d}\left( \mathbb{H}\right) \right) =\{0\}$ if $r$
is not a multiple of $4$. In $\mathrm{Gr}_{d}\left( \mathbb{H}\right) $
there is just one $4$-dimensional cell which is the orbit $\mathrm{Sl}\left(
2,\mathbb{H}\right) _{d,d+1}V_{d}$. Here,  $\mathrm{Sl}\left( 2,\mathbb{H}%
\right) _{d,d+1}=\langle \exp \mathfrak{sl}\left( 2,\mathbb{H}\right)
_{d,d+1}\rangle \approx \mathrm{Sl}\left( 2,\mathbb{H}\right) $ and $%
\mathfrak{sl}\left( 2,\mathbb{H}\right) _{d,d+1}$ is the algebra of
matrices with nonzero entries only in the entries $\left( d,d\right) $, $%
\left( d,d+1\right) $, $\left( d+1,d\right) $ and $\left( d+1,d+1\right) $.
Analogous to the Lemma \ref{lemEsse4} we have that $\mathrm{Sl}\left( 2,%
\mathbb{H}\right) _{d,d+1}V_{d}$ is diffeomorphic to $S^{4}$.

Now, by the Hurewicz homomorphism $\pi _{4}\left( \mathrm{Gr}_{d}\left( 
\mathbb{H}\right) \right) \approx H_{4}\left( \mathrm{Gr}_{d}\left( \mathbb{H%
}\right) \right) $ because the homology is trivial in degrees less than $4$.
It follows that $\pi _{4}\left( \mathrm{Gr}_{d}\left( \mathbb{H}\right)
\right)$ $\approx$ $\mathbb{Z}$ and the equivalence class of the orbit $\mathrm{%
Sl}\left( 2,\mathbb{H}\right) _{d,d+1}V_{d}\approx S^{4}$ is a generator of $%
\pi _{4}\left( \mathrm{Gr}_{d}\left( \mathbb{H}\right) \right) $. The next
lemma shows that $\mathrm{Sl}\left( 2,\mathbb{H}\right) _{1,n}V_{d}\approx
S^{4}$ is a generator as well.

\begin{lema}
The orbits $\mathrm{Sl}\left( 2,\mathbb{H}\right) _{d,d+1}V_{d}\approx S^{4}$
and $\mathrm{Sl}\left( 2,\mathbb{H}\right) _{1,n}V_{d}\approx S^{4}$ are
homotopic to each other.
\end{lema}

\begin{profe}
The homotopy is performed by the product of two one-parameter subgroups. Let 
$A,B\in \mathfrak{sl}\left( n,\mathbb{H}\right) $ be the matrices such that $%
Ae_{1}=e_{d}$, $Ae_{d}=-e_{1}$, $Be_{d+1}=e_{n}$, $Be_{n}=-e_{d+1}$ and $%
Ae_{r}=Be_{r}=0$ elsewhere. Put $P\left( t\right) =e^{tA}e^{tB}$. Then for
all $t$, $P\left( t\right) V_{d}=V_{d}$ and $P\left( \pi /2\right) $
permutes the subspaces spanned by $\{e_{d},e_{d+1}\}$ and $\{e_{1},e_{n}\}$
so that $P\left( \pi /2\right) \mathrm{Sl}\left( 2,\mathbb{H}\right)
_{d,d+1}P\left( \pi /2\right) ^{-1}=\mathrm{Sl}\left( 2,\mathbb{H}\right)
_{1,n}$. Hence
\begin{eqnarray*}
	P\left( \pi /2\right) \mathrm{Sl}\left( 2,\mathbb{H}\right)
	_{d,d+1}V_{d}&=&P\left( \pi /2\right) \mathrm{Sl}\left( 2,\mathbb{H}\right)
	_{d,d+1}P\left( \pi /2\right) ^{-1}P\left( \pi /2\right) V_{d}\\ &=&\mathrm{Sl}%
	\left( 2,\mathbb{H}\right) _{1,n}V_{d} ,
\end{eqnarray*}
showing that the map $t\mapsto P\left( t\right) \mathrm{Sl}\left( 2,\mathbb{H%
}\right) _{d,d+1}V_{d}$ is a homotopy between the orbits $\mathrm{Sl}\left(
2,\mathbb{H}\right) _{d,d+1}V_{d}$ and $\mathrm{Sl}\left( 2,\mathbb{H}%
\right) _{1,n}V_{d}$.
\end{profe}

Now we can start the proof of Theorem \ref{teoTopology} in the case when $%
\mathrm{Sl}\left( 2,\mathbb{H}\right) _{1,n}\subset S$. Denote by $N$ the
nilpotent group of lower triangular matrices in $\mathrm{Sl}\left( n,\mathbb{%
H}\right) $ having $1$'s at the diagonal. It is well known (and easy to
prove) that $NV_{d}$ is an open and dense set in $\mathrm{Gr}_{d}\left( 
\mathbb{H}\right) $. Hence $NV_{d}\cap C_{d}\neq \emptyset $ because $%
\mathrm{int}C_{d}\neq \emptyset $ where $C_{d}$ is the invariant control set
of $S$ in $\mathrm{Gr}_{d}\left( \mathbb{H}\right) $.

The assumption $\mathrm{Sl}\left( 2,%
\mathbb{H}\right) _{1,n}\subset S$ of Theorem \ref{teoTopology}  implies that $gC_{d}\subset C_{d}$ for
any $g\in \mathrm{Sl}\left( 2,\mathbb{H}\right) _{1,n}$. Since $C_{d}$ is
closed it follows that any limit $\lim g_{l}x$ with $x\in C_{d}$ and $%
g_{l}\in \mathrm{Sl}\left( 2,\mathbb{H}\right) _{1,n}$ also belongs to $C_{d}$.

Now, take $x=gV_{d}\in NV_{d}\cap C_{d}$ with $g\in N$ and 
\begin{equation*}
h=\mathrm{diag}\{\lambda ,1,\ldots ,1,\lambda ^{-1}\}\in \mathrm{Sl}\left( 2,%
\mathbb{H}\right) _{1,n}
\end{equation*}%
with $\lambda >1$. As $l\rightarrow +\infty $ the sequence of conjugations $%
h^{l}gh^{-l}$ converges to the matrix $g_{1}\in N$ that has zeros at the
first column and the last row outside the diagonal. We have $h^{-l}V_{d}=V_{d}$
so that $h^{l}x=h^{l}gV_{d}=h^{l}gh^{-l}V_{d}$ implying that $W=\lim
h^{l}V_{d}=g_{1}V_{d}\in C_{d}$. Therefore the orbit $\mathrm{Sl}\left( 2,%
\mathbb{H}\right) _{1,n}W$ is entirely contained in $C_{d}$.

The next step is to prove that the orbit $\mathrm{Sl}\left( 2,\mathbb{H}%
\right) _{1,n}W\subset C_{d}$ is a sphere $S^{4}$ homotopic to $\mathrm{Sl}%
\left( 2,\mathbb{H}\right) _{1,n}V_{d}$. By the zeros in the first column
and the last row of $g_{1}$, the subspace $W=g_{1}V_{d}$ is the direct sum $%
\langle e_{1}\rangle \oplus W_{1}$ where $W_{1}$ is a $\left( d-1\right) $%
-dimensional subspace of $\mathrm{span}_{\mathbb{H}}\{e_{2},\ldots ,e_{n-1}\}
$. If $d=1$ then $W=V_{d}$ and we are done. Otherwise, let $G=\mathrm{Sl}%
\left( n-2,\mathbb{H}\right) _{2,\ldots ,n-1}\approx \mathrm{Sl}\left( n-2,%
\mathbb{H}\right) $ be the subgroup of matrices in $\mathrm{Sl}\left( n,\mathbb{H%
}\right) $ whose restriction to $\mathrm{span}_{\mathbb{H}}\{e_{1},e_{n}\}$
is the identity. Then $G$ commutes with $\mathrm{Sl}\left( 2,\mathbb{H}%
\right) _{1,n}$ so that if $g\in G$ then $g\mathrm{Sl}\left( 2,\mathbb{H}%
\right) _{1,n}W=\mathrm{Sl}\left( 2,\mathbb{H}\right) _{1,n}gW$, that is,
the image under $g\in G$ of the orbit $\mathrm{Sl}\left( 2,\mathbb{H}\right)
_{1,n}W$ is again an orbit of $\mathrm{Sl}\left( 2,\mathbb{H}\right) _{1,n}$%
. Moreover, $G$ acts transitively in the Grasmmannian of $\left( d-1\right) $%
-dimensional subspaces of $\mathrm{span}_{\mathbb{H}}\{e_{2},\ldots
,e_{n-1}\}$. Hence there exists $g\in G$ such that $gW_{1}=\mathrm{span}_{%
\mathbb{H}}\{e_{2},\ldots ,e_{d}\}$ so that $gW=V_{d}$. It follows that the
orbit $\mathrm{Sl}\left( 2,\mathbb{H}\right) _{1,n}W$ is diffeomorphic to $%
\mathrm{Sl}\left( 2,\mathbb{H}\right) _{1,n}V_{d}$ and hence is a sphere $%
S^{4}$. Furthermore, $G$ is connected so that there is a continuous curve $%
g_{t}\in G$ with $g_{0}=1$ and $g_{1}=g$. Hence $t\mapsto g_{t}\mathrm{Sl}%
\left( 2,\mathbb{H}\right) _{1,n}W$ is a homotopy between $\mathrm{Sl}\left(
2,\mathbb{H}\right) _{1,n}W$ and $\mathrm{Sl}\left( 2,\mathbb{H}\right)
_{1,n}V_{d}$ showing that $\mathrm{Sl}\left( 2,\mathbb{H}\right)
_{1,n}W\approx S^{4}$ is not contractible.

We proved that the invariant control set $C_{d}$ of the semigroup $S$ in $%
\mathrm{Gr}_{d}\left( \mathbb{H}\right) $ contains a noncontractible sphere $%
S^{4}$. Hence $C_{d}$ is not contractible either. Since $d=1,\ldots ,n-1$ is
arbitrary, $S$ cannot be a proper semigroup, concluding the proof of Theorem %
\ref{teoTopology} in the case when $\mathrm{Sl}\left( 2,\mathbb{H}\right)
_{1,n}\subset S$.

\begin{coro}
\label{corqualqueraiz}Let $T\subset \mathrm{Sl}\left( n,\mathbb{H}\right) $
be a semigroup with nonempty interior. Suppose that $\mathrm{Sl}\left( 2,%
\mathbb{H}\right) _{r,s}\subset T$ for some pair $\left( r,s\right) $, $%
r\neq s$. Then $T=\mathrm{Sl}\left( n,\mathbb{H}\right) $.
\end{coro}

\begin{profe}
Let $P$ be a matrix in $\mathrm{Sl}\left( n,\mathbb{H}\right) $ that permutes
the subspaces $\langle e_{1}\rangle $ and $\langle e_{r}\rangle $ and the
subspaces $\langle e_{n}\rangle $ and $\langle e_{s}\rangle $. Then $P%
\mathrm{Sl}\left( 2,\mathbb{H}\right) _{r,s}P^{-1}=\mathrm{Sl}\left( 2,%
\mathbb{H}\right) _{1,n}$ so that the semigroup $PTP^{-1}$ contains $\mathrm{%
Sl}\left( 2,\mathbb{H}\right) _{1,n}$. Since $\mathrm{int}PTP^{-1}\neq
\emptyset $ we conclude, by Theorem \ref{teoTopology}, that $PTP^{-1}=\mathrm{%
Sl}\left( n,\mathbb{H}\right) $, hence $T=\mathrm{Sl}\left( n,\mathbb{H}%
\right) $.
\end{profe}

\section{Application to Controllability}

In this section we apply Theorem \ref{teoTopology} to show the following sufficient conditions
 for controllability of an invariant system like (\ref{forsis}) in Sl$(n,\mathbb{H})$. The statement of this theorem is, in some
 sense, a  $\mathfrak{sl}\left( n,\mathbb{H}\right)$-version of  a well-known approach to controllability (see e.g., El Alssoudi, Gauthier and Kupka \cite{EGK}).

\begin{teo}
	\label{teoregularreal}Let $\dot{g}=A\left( g\right) +uB\left( g\right) $ be
	an invariant control system in $\mathrm{Sl}\left( n,\mathbb{H}\right) $
	where $A,B\in \mathfrak{sl}\left( n,\mathbb{H}\right) $. Such a system with
	unrestricted controls ($u\in \mathbb{R}$) is controllable if the following
	conditions are satisfied.
	
	\begin{enumerate}
		\item[H1.] The pair $\left( A,B\right) $ generates $\mathfrak{sl}\left( n,\mathbb{%
			H}\right) $ as a Lie algebra (Lie algebra rank condition).
		
		\item[H2.] $B=\mathrm{diag}\{a_{1}+ib_{1},\ldots ,a_{n}+ib_{n}\}$ with $%
		a_{1}>a_{2}\geq \cdots \geq a_{n-1}>a_{n}$, $b_{n}\neq 0\neq b_{1}$ and $%
		b_{1}/b_{n}$ is irrational.
		
		
		\item[H3.] Denote the $1,n$ and $n,1$ entries of the matrix $A$ by $p\in \mathbb{H%
		}$ and $q\in \mathbb{H}$, respectively. Let $\mathbb{H}_{1,i}$ and $\mathbb{H%
		}_{j,k}$ be the (real) subespaces of $\mathbb{H}$ spanned by $\{1,i\}$ and $%
		\{j,k\}$ respectively. Then $p$ and $q$ do not belong to $\mathbb{H}%
		_{1,i}\cup \mathbb{H}_{j,k}$.
	\end{enumerate}
\end{teo}

The proof of the Theorem \ref{teoregularreal} will be made throughout this section and it is immediate from Theorem \ref{teoTopology} combined with the following proposition ensuring that $\mathrm{Sl}\left( 2,%
\mathbb{H}\right) _{1,n}$ is contained in the control semigroup of the system. Although the Lie algebra rank condition will not be needed for this proposition, it allows us to conclude the proof of the Theorem \ref{teoregularreal} by ensuring that the control semigroup $S$ has nonempty interior, leading us to the conditions required for Theorem \ref{teoTopology}.

\begin{prop}
	\label{teoHighestRoot}Under the conditions H2 and H3 of Theorem \ref{teoregularreal},
	the semigroup $S$ of the system contains the group $\mathrm{Sl}\left( 2,%
	\mathbb{H}\right) _{1,n}$.
\end{prop}

To prove this proposition, let $S$ be the control semigroup for the invariant system (\ref%
{forsis}) and write 
\begin{equation*}
\mathfrak{c}\left( S\right) =\{X\in \mathfrak{sl}\left( n,\mathbb{H}\right)
:\forall t\geq 0,~e^{tX}\in \mathrm{cl}S\}
\end{equation*}%
for the Lie wedge of $S$ (see \cite{JK}, \cite{JK1} and Hilgert, Hofmann and Lawson \cite{hhl}). The
main properties of $\mathfrak{c}\left( S\right) $ are: 

1) $\mathfrak{c}%
\left( S\right) $ is a closed convex cone in the Lie algebra $\mathfrak{sl}%
\left( n,\mathbb{H}\right) $; 

2) $\mathfrak{c}\left( S\right) \cap \left( -%
\mathfrak{c}\left( S\right) \right) $ is a Lie subalgebra and 

3) If $X\in 
\mathfrak{c}\left( S\right) \cap \left( -\mathfrak{c}\left( S\right) \right) 
$ then $e^{\mathrm{ad}\left( X\right) }\mathfrak{c}\left( S\right) =%
\mathfrak{c}\left( S\right) $.

By definition of $S$ we have that $A+uB\in \mathfrak{c}\left( S\right) $ for
all $u\in \mathbb{R}$ (since we consider unrestricted controls). Hence $A\in 
\mathfrak{c}\left( S\right) $ and if $u\neq 0$ then 
\begin{equation*}
\frac{1}{\left\vert u\right\vert }A+\frac{u}{\left\vert u\right\vert }B=%
\frac{1}{\left\vert u\right\vert }\left( A+uB\right) \in \mathfrak{c}\left(
S\right) .
\end{equation*}%
Taking limits as $u\rightarrow \pm \infty $ we see that $\pm B\in \mathfrak{c%
}\left( S\right) $, that is, $B\in \mathfrak{c}\left( S\right) \cap \left( -%
\mathfrak{c}\left( S\right) \right) $. It follows that $e^{t\mathrm{ad}%
	\left( B\right) }A\in \mathfrak{c}\left( S\right) $ and hence $e^{-t\left(
	a_{1}-a_{n}\right) }e^{t\mathrm{ad}\left( B\right) }A\in \mathfrak{c}\left(
S\right) $ for all $t\in \mathbb{R}$ where $a_{1},\ldots ,a_{n}$ are the
real parts of the entries of $B$.

Now by assumption we have $a_{1}>a_{2}>\cdots >a_{n}$ so that as $t\rightarrow
+\infty $ the entries $e^{-t\left( a_{1}-a_{n}\right) }e^{t\mathrm{ad}\left(
	B\right) }A$ converge to $0$ except for the $\left( 1,n\right) $-entry. The $%
\left( 1,n\right) $-entry of $e^{-t\left( a_{1}-a_{n}\right) }e^{t\mathrm{ad}%
	\left( B\right) }A$ is $e^{it\left( b_{1}-b_{n}\right) }p$ where $p$ is as
in the statement of the theorem and $b_{1},\ldots ,b_{n}$ are the imaginary
parts of the entries of $B$. Choosing a sequence $t_{k}\rightarrow +\infty $
such that $e^{it\left( b_{1}-b_{n}\right) }\rightarrow 1$ we conclude that 
\begin{equation*}
X=\left( 
\begin{array}{ccc}
0 & \cdots & p \\ 
\vdots & \ddots & \vdots \\ 
0 & \cdots & 0%
\end{array}%
\right) \in \mathfrak{c}\left( S\right) .
\end{equation*}%
Using again the properties of $\mathfrak{c}\left( S\right) $ as a Lie wedge
we have that for all $t,s\in \mathbb{R}$, 
\begin{equation}
e^{-t\left( a_{1}-a_{n}\right) }e^{t\mathrm{ad}\left( B\right) }X=\left( 
\begin{array}{ccc}
0 & \cdots & e^{itb_{1}}pe^{-itb_{n}} \\ 
\vdots & \ddots & \vdots \\ 
0 & \cdots & 0%
\end{array}%
\right) \in \mathfrak{c}\left( S\right) .  \label{forconjugaupper}
\end{equation}
The following lemma about conjugation of quaternions shows that $\mathfrak{g}%
_{\alpha _{1n}}=\mathrm{span}_{\mathbb{H}}\{X\}$ is contained in $\mathfrak{c%
}\left( S\right) $.

\begin{lema}
	Consider the action of the circle group $S^{1}=\{e^{it}:t\in \left[ 0,2\pi %
	\right] \}$ in $\mathbb{H}$ given by conjugation $\left( t,q\right) \mapsto
	e^{it}qe^{-it}$. Write $q=a+b$ with $a=x_{1}+ix_{2}\in \mathbb{H}_{\{1,i\}}$
	and $b=jx_{3}+kx_{4}\in \mathbb{H}_{\{j,k\}}$. Suppose that $a\neq 0\neq b$,
	that is, $q\notin \mathbb{H}_{\{1,i\}}\cup \mathbb{H}_{\{j,k\}}$. Then the
	orbit $T^{2}q$ is a $2$-dimensional torus and $\mathbb{H}$ is the convex
	cone generated by $T^{2}q$.
\end{lema}

\begin{profe}
	For the first statement it suffices to prove that the orbit is $2$%
	-dimensional because it is a quotient of $T^{2}$. For this purpose we note
	that the tangent space of $T^{2}q$ at $q$ is spanned by 
	\begin{equation*}
	\frac{\partial }{\partial t}\left( e^{it}qe^{-is}\right) _{\left\vert \left(
		0,0\right) \right. }=iq\qquad \frac{\partial }{\partial s}\left(
	e^{it}qe^{-is}\right) _{\left\vert \left( 0,0\right) \right. }=-qi.
	\end{equation*}
	Now, 
	\begin{eqnarray*}
		iq &=&ix_{1}-x_{2}+kx_{3}-jx_{4}=v+w \\
		-qi &=&-ix_{1}+x_{2}+kx_{3}-jx_{4}=-v+w
	\end{eqnarray*}%
	with $v=ix_{1}-x_{2}$ and $w=kx_{3}-jx_{4}$. The assumption about $q$ says
	that $v\neq 0\neq w$ so that $\{v,w\}$ is  linearly independent. Hence $\{iq,qi\}$ is 
	 linearly independent as well because $2v=iq+qi$ and $2w=iq-qi$.
	
	To the convex cone $C$ generated by $T^{2}q$ take $r=e^{it}qe^{-is}\in
	T^{2}q$. Then $-r=e^{i\pi }r=e^{i\left( t+\pi \right) }qe^{-is}$ also
	belongs to $T^{2}q$. Hence $C$ is a subspace. The orbit contains $e^{i\pi
		/2}q=iq=ia+ib$ and $qe^{-i3\pi /2}=qi=ai+bi=ia-ib$. So that $C$ contains $ia$
	and $ib$ and hence contains $\mathbb{H}_{\{1,i\}}$ and $\mathbb{H}_{\{j,k\}}$
	because $C$ is invariant by left multiplication by $i$. Thus $\mathbb{H}$ is
	the cone gerated by $T^{2}q$.
\end{profe}

\begin{coro}
	Let $c_{1},c_{2}\in \mathbb{R}$ with $c_{1}c_{2}\neq 0$ and $c_{1}/c_{2}$
	irrational. Take $q\in \mathbb{H}$ with $q\notin \mathbb{H}_{\{1,i\}}\cup 
	\mathbb{H}_{\{j,k\}}$. Then $\mathbb{H}$ is the closed convex cone generated
	by the curve $e^{itc_{1}}qe^{-itc_{2}}$.
\end{coro}

\begin{profe}
	Since $c_{1}/c_{2}$ is irrational, the curve $\left(
	e^{itc_{1}},e^{itc_{2}}\right) $ is dense in the torus $T^{2}$. Hence $%
	t\mapsto e^{itc_{1}}qe^{-itc_{2}}$ is dense in the orbit $T^{2}q$ which
	implies the corollary.
\end{profe}

Applying this corollary to the curve (\ref{forconjugaupper}) it follows, by
the assumption on $B$ in Theorem \ref{teoregularreal}, that the subspace $%
\mathfrak{g}_{\alpha _{1n}}=\mathrm{span}_{\mathbb{H}}\{X\}$ is contained in 
$\mathfrak{c}\left( S\right) $ and hence in $\mathfrak{c}\left( S\right)
\cap \left( -\mathfrak{c}\left( S\right) \right) $.

By similar arguments we get lower triangular matrices in $\mathfrak{c}%
\left( S\right) $: taking limits as $t\rightarrow -\infty $ of $e^{-t\left(
	a_{1}-a_{n}\right) }e^{t\mathrm{ad}\left( B\right) }A$ it follows that 
\begin{equation*}
Y=\left( 
\begin{array}{ccc}
0 & \cdots & 0 \\ 
\vdots & \ddots & \vdots \\ 
q & \cdots & 0%
\end{array}%
\right) \in \mathfrak{c}\left( S\right) ,
\end{equation*}%
hence applying the same idea  we conclude that $\mathfrak{g%
}_{\alpha _{n1}}=\mathrm{span}_{\mathbb{H}}\{Y\}$ is contained in $\mathfrak{%
	c}\left( S\right) $ and hence in $\mathfrak{c}\left( S\right) \cap \left( -%
\mathfrak{c}\left( S\right) \right) $.

Now, the Lie algebra generated by $\mathfrak{g}_{\alpha _{1n}}$ and $%
\mathfrak{g}_{\alpha _{n1}}$ is $\mathfrak{sl}\left( 2,\mathbb{H}\right)
_{1n}$ so that this Lie algebra is contained in $\mathfrak{c}\left( S\right) 
$. It follows that $\mathrm{Sl}\left( 2,\mathbb{H}\right) _{1n}$ is
contained in $S$, concluding the proof of Proposition \ref{teoHighestRoot}.

\section{Cartan subalgebras and genericity\label{secgeneric}}

In this section we prove that controllability for invariant control systems in ${\rm Sl}\left(n,\mathbb{H}\right)$ is a generic property. To prove it we first show that the set of pairs conjugate to a pair $\left( A,B\right) $ satisfying the conditions of Theorem \ref{teoregularreal} is dense in $\mathfrak{sl}\left( n,\mathbb{H}\right) ^{2}$.

We start by observing that the algebra of diagonal matrices 
\begin{equation}
\mathfrak{h}=\{\mathrm{diag}\{a_{1}+ib_{1},\ldots
,a_{n}+ib_{n}\}:a_{r},b_{r}\in \mathbb{R},~a_{1}+\cdots +a_{n}=0\}
\label{forCartansubalg}
\end{equation}%
is a Cartan subalgebra of $\mathfrak{sl}\left( n,\mathbb{H}\right) $ since
it is maximal abelian and $\mathrm{ad}\left( H\right) $ is semi-simple for
any $H\in \mathfrak{h}$.

Next we prove  that up to conjugation, $\mathfrak{h}$ is the only Cartan
subalgebra of $\mathfrak{sl}\left( n,\mathbb{H}\right) $. Then we
recall the Cartan decomposition $\mathfrak{sl}\left( n,\mathbb{H}\right) =%
\mathfrak{sp}\left( n\right) \oplus \mathfrak{s}$ where $\mathfrak{s}$ is
the subspace of Hermitian quaternionic matrices in $\mathfrak{sl}\left( n,%
\mathbb{H}\right) $. The subspace $\mathfrak{a}\subset \mathfrak{s}$ of real
diagonal matrices with zero trace is a maximal abelian subalgebra contained
in $\mathfrak{s}$.

Note that the Cartan subalgebra $\mathfrak{h}$ decomposes as $\mathfrak{h}%
=\left( \mathfrak{h\cap sp}\left( n\right) \right) \oplus \mathfrak{a}$.
More generally $\mathfrak{j}$ is said to be a standard Cartan subalgebra if
it decomposes as $\mathfrak{j}=\mathfrak{j}_{\mathfrak{k}}\oplus \mathfrak{j}%
_{\mathfrak{a}}$ with $\mathfrak{j}_{\mathfrak{k}}=\mathfrak{j}\cap 
\mathfrak{k}$ and $\mathfrak{j}_{\mathfrak{a}}=\mathfrak{j}\cap \mathfrak{a}$%
. The following statement is a basic fact for the classification of Cartan
subalgebras in real semi-simple Lie algebras (Theorem of Kostant-Sugiura).

\begin{prop}
	Any Cartan subalgebra of $\mathfrak{g}=\mathfrak{sl}\left( n,\mathbb{H}%
	\right) $ is conjugate (by an inner automorphism) to a standard Cartan
	subalgebra $\mathfrak{j}$.
\end{prop}

\begin{profe}
	See Warner \cite{W}, Section 1.3.1.
\end{profe}

In the next proposition, we prove that in $\mathfrak{sl}\left( n,\mathbb{H}%
\right)$ there is a unique conjugacy class of Cartan subalgebras. We give a
direct proof without relying in the general classification theorem (Theorem
of Kostant-Sugiura).

\begin{prop}
	\label{propsubalgCartan}Every Cartan subalgebra of $\mathfrak{sl}\left( n,%
	\mathbb{H}\right) $ is conjugate (by an inner automorphism) to the
	subalgebra $\mathfrak{h}$ defined in (\ref{forCartansubalg}).
\end{prop}

\begin{profe}
	Let $\mathfrak{j}=\mathfrak{j}_{\mathfrak{k}}\oplus \mathfrak{j}_{\mathfrak{a%
	}}$ be a standard Cartan subalgebra. The following simple arguments show
 that $\mathfrak{j}_{\mathfrak{k}}$ is a Cartan subalgebra of 
	$\mathfrak{sp}\left( n\right) $ and $\mathfrak{j}_{\mathfrak{a}}=\mathfrak{a}
	$. We have $\dim \mathfrak{j}_{\mathfrak{a}}\leq \dim \mathfrak{a}=n-1$.
	Also, $\dim \mathfrak{j}_{\mathfrak{k}}\leq \mathrm{rank} \mathfrak{sp}\left(
	n\right) =n$ because $\mathfrak{j}_{\mathfrak{k}}$ is an abelian subalgebra
	of $\mathfrak{sp}\left( n\right) $ and hence is contained in a Cartan
	subalgebra of $\mathfrak{sp}\left( n\right) $ whose dimension is $\mathrm{%
		rank} \mathfrak{sp}\left( n\right) $. On the other hand $\dim \mathfrak{j}=%
	\mathrm{rank} \mathfrak{sl}\left( n,\mathbb{H}\right) =\dim \mathfrak{h}=2n-1$%
	. Hence we must have $\dim \mathfrak{j}_{\mathfrak{k}}=n$ and $\dim 
	\mathfrak{j}_{\mathfrak{a}}=n-1$. By the first equality $\mathfrak{j}_{%
		\mathfrak{k}}$ is a Cartan subalgebra of $\mathfrak{sp}\left( n\right) $
	while the second equality shows that $\mathfrak{j}_{\mathfrak{a}}=\mathfrak{a%
	}$.
	
	Now $\mathfrak{j}_{\mathfrak{k}}$ commutes with $\mathfrak{a}$ and hence is
	contained in the algebra $\mathfrak{m}\approx \mathfrak{sp}\left( 1\right)
	^{n}$ of diagonal matrices with entries in the imaginary quaternions $
	\mathrm{Im}\mathbb{H}$. Since $\dim \mathfrak{j}_{\mathfrak{k}}=n=%
	\mathrm{rank}\mathfrak{sp}\left( 1\right) ^{n}$ it follows that there is an
	inner automorphism $g=e^{\mathrm{ad}\left( X\right) }$, $X\in \mathfrak{m}$,
	such that $g\left( \mathfrak{j}_{\mathfrak{k}}\right) =\mathfrak{h}_{k}$ and 
	$g$ fixes $\mathfrak{a}$. Therefore $g\left( \mathfrak{j}\right) =\mathfrak{h%
	}$ showing that any standard Cartan subalgebra is conjugate to $\mathfrak{h}$%
	. By the above proposition, $\mathfrak{h}$ is a representative of the unique
	conjugacy class of Cartan subalgebras of $\mathfrak{sl}\left( n,\mathbb{H}%
	\right) $.
\end{profe}

 Now denote by $\mathfrak{a}^{+}$ the
Weyl chamber of real diagonal matrices 
\begin{equation*}
\mathrm{diag}\{a_{1},\ldots ,a_{n}\}\qquad a_{1}>\cdots >a_{n}.
\end{equation*}%
A matrix $B$ satifying the second condition of Theorem \ref{teoregularreal}
belongs to $\mathfrak{a}^{+}+\mathfrak{h}_{\mathfrak{k}}$ where $\mathfrak{h}%
_{\mathfrak{k}}$ is as above the space of diagonal matrices with entries in $%
i\mathbb{R}$. Denote by $D_{0}\subset \mathfrak{a}^{+}+\mathfrak{h}_{%
	\mathfrak{k}}$ the set of the matrices $B$ satisfying that condition. By
definition if $H\in \mathfrak{a}^{+}$ and $X=\mathrm{diag}\{ib_{1},\ldots
,ib_{n}\}\in \mathfrak{h}_{\mathfrak{k}}$ then $H+X\in D_{0}$ if and only if 
$b_{1}b_{n}\neq 0$ and $b_{1}/b_{n}$ is irrational. Hence $D_{0}$ is a dense
subset of $\mathfrak{a}^{+}+\mathfrak{h}_{\mathfrak{k}}$.

Let $\mathcal{W}$ be the permutation group in $n$ letters (Weyl group)
acting on the diagonal matrices by permutation of indices. The set of
translates $\mathcal{W}\mathfrak{a}^{+}=\{w\mathfrak{a}^{+}:w\in \mathcal{W}%
\}$ is open and dense in $\mathfrak{a}$. Since $D_{0}$ is dense in $%
\mathfrak{a}^{+}+\mathfrak{h}_{\mathfrak{k}}$ it follows that 
\begin{equation*}
\mathcal{W}D_{0}=\{wD_{0}:w\in \mathcal{W}\}\subset \mathfrak{a}+\mathfrak{h}%
_{\mathfrak{k}}=\mathfrak{h}
\end{equation*}%
is dense in $\mathfrak{h}$.

We apply now Proposition \ref{propsubalgCartan} ensuring that every Cartan
subalgebra is conjugate to $\mathfrak{h}$. This implies that the set $\{%
\mathrm{Ad}\left( g\right) \mathfrak{h}:g\in \mathrm{Sl}\left( n,\mathbb{H}%
\right) \}$ is dense in $\mathfrak{sl}\left( n,\mathbb{H}\right) $ because
the set of regular elements is dense and each regular element is contained
in a Cartan subalgebra. With these facts in mind we get the following
density result.

\begin{prop}
	Let $D$ be the set of conjugates of the matrices $B$ satisfying the second
	condition of Theorem \ref{teoregularreal}. Then $D$ is dense in $\mathfrak{sl%
	}\left( n,\mathbb{H}\right) $.
	\label{denseD}
\end{prop}

\begin{profe}
	Take an open set $U\subset \mathfrak{sl}\left( n,\mathbb{H}\right) $. Then
	there exists a regular element $X$ of $\mathfrak{sl}\left( n,\mathbb{H}%
	\right) $ with $X\in U$. Let $\mathfrak{h}_{X}$ be the unique Cartan
	subalgebra containing $X$. By Proposition \ref{propsubalgCartan} there
	exists $g\in \mathrm{Sl}\left( n,\mathbb{H}\right) $ such that $\mathrm{Ad}%
	\left( g\right) \mathfrak{h}_{X}=\mathfrak{h}$. So that $\mathrm{Ad}\left(
	g\right) U\cap \mathfrak{h}$ is a nonempty open set of $\mathfrak{h}$ and
	hence $\mathrm{Ad}\left( g\right) U\cap \mathcal{W}D_{0}\neq \emptyset $.
	This means that $U$ meets $\mathrm{Ad}\left( g^{-1}\right) \mathcal{W}%
	D_{0}\subset D$. Since $U$ is arbitrary this proves that $D$ is dense.
\end{profe}

Now we can show the main result of this section.

\begin{teo}
	\label{teogeneric}
	There is an open and dense set $C\subset \mathfrak{sl}%
	\left( n,\mathbb{H}\right) ^{2}$ such that the control system $\dot{g}%
	=A\left( g\right) +uB\left( g\right) $ with unrestricted controls ($u\in 
	\mathbb{R}$) is controllable for all pairs $\left( A,B\right) \in C$.
\end{teo}

To prove this theorem, first note that the union $\mathbb{H}_{1,i}\cup\mathbb{H}_{j,k}$ is a nowhere dense subset of $\mathbb{H}$, which implies that its complement is an open and dense subset of $\mathbb{H}$. Consequently, the set of matrices $A$ satisfying the third condition of Theorem \ref{teoregularreal} is open and dense in $\mathfrak{sl}\left( n,\mathbb{H}\right)$.

\begin{rem}
	For $\Omega\subset M\times N$ open, the set $\pi_1\left(\Omega\cap\pi_2^{-1}(b)\right)\subset M$ is open in $M$ for any $b\in N$. Here, $M$ and $N$ are arbitrary metric spaces and $\pi_1:M\times N\to N$ and $\pi_2:M\times N\to N$ are the canonical projections in the first and second coordinates, respectively. To see this just let the continuous map $i_b:M\to M\times N$, $i_b(x)=(x,b)$, and observe that
	$$\pi_1\left(\Omega\cap\pi_2^{-1}(b)\right)=\left\{\pi_1(x,b)\;|\;(x,b)\in\Omega\right\}=\left\{x\in M\;|\;i_b(x)\in\Omega\right\}=(i_b)^{-1}(\Omega).$$
\end{rem}

We can now prove that the set $C\subset \mathfrak{sl}(n,\mathbb{H})^2$ of the conjugates of pairs satisfying the three conditions of the Theorem \ref{teoregularreal} is dense in $\mathfrak{sl}(n,\mathbb{H})^2$.

So, let $O$ be an open subset of $\mathfrak{sl}(n,\mathbb{H})^2$. Since the set of pairs $(A,B)$ satisfying H1 is open and dense in $\mathfrak{sl}(n,\mathbb{H})^2$, there is $(A,B)\in O$ satisfying H1. Further, there exists $O'\ni(A,B)$ for which every pair belonging to $O'$ satisfies H1. Without loss of generality we can assume $O'\subset O$. Now, $\pi_2(O')$ is open in $\mathfrak{sl}(n,\mathbb{H})$ and by the Proposition \ref{denseD} we can choose $\tilde{B}\in\pi_2(O')\cap D$. As the set $\pi_1\left(O'\cap\pi_2^{-1}(\tilde{B})\right)$ is open in $\mathfrak{sl}(n,\mathbb{H})^2$, by the above considerations we can take $\tilde{A}\in\pi_1\left(\pi_2^{-1}(\tilde{B})\cap O'\right)$ satisfying H3. Thus the pair $(\tilde{A},\tilde{B})$ has the following properties:
\begin{itemize}
	\item[i)] $(\tilde{A},\tilde{B})\in O'\subset O$.
	\item[ii)] $(\tilde{A},\tilde{B})$ is conjugate to a pair satisfying H1, H2 and H3.
\end{itemize}
That is, $(\tilde{A},\tilde{B})\in O\cap C$ proving that $C$ is dense in $\mathfrak{sl}(n,\mathbb{H})^2$. Finally, as invariant systems remain controllable under small perturbations we can slightly enlarge the dense set $C$ to get the open and dense set as claimed in Theorem \ref{teogeneric}.

\end{document}